\documentclass[10pt,english]{amsart}
\usepackage{helvet}
\usepackage[T1]{fontenc}
\usepackage[latin1]{inputenc}
\usepackage{geometry}
\usepackage{pifont}
\usepackage{amssymb}

\numberwithin{equation}{section} 
\numberwithin{figure}{section} 

\usepackage{latexsym,amsfonts, xypic,amssymb,amsthm}
\xyoption{all}
  \theoremstyle{plain}
  \newtheorem{thm}{Theorem}[section]
  \theoremstyle{plain}
  \newtheorem*{thm*}{Theorem}
  \theoremstyle{remark}
  \newtheorem*{acknowledgement*}{Acknowledgement}
  \theoremstyle{remark}
  \newtheorem{rem}[thm]{Remark}
  \theoremstyle{remark}
  
  \theoremstyle{definition}
  \newtheorem{enavant}[thm]{}

\usepackage{babel}

\begin{document}

\title{The cylinder over the Koras-Russell cubic threefold has a trivial Makar-Limanov invariant }

\author{Adrien Dubouloz}

\address{Cnrs, Institut de Math\'ematiques de Bourgogne, Universit\'e de
Bourgogne, 9 avenue Alain Savary - BP 47870, 21078 Dijon cedex, France}

\email{adrien.dubouloz@u-bourgogne.fr}


\begin{abstract} We show that  the cylinder $X\times\mathbb{A}^1$ over the  Koras-Russell cubic threefold  $$X=\{x+x^2y+z^2+t^3=0\}\subset \mathbb{A}^4$$ has a trivial Makar-Limanov invariant ${\rm ML}(X\times \mathbb{A}^1)=\mathbb{C}$. This means equivalently that the only regular functions on $X\times\mathbb{A}^1$ that are invariant under \emph{all} algebraic actions of the additive  group $\mathbb{G}_a$ on $X\times\mathbb{A}^1$ are constants. \\ \end{abstract}
\maketitle

\section*{Introduction}

The \emph{Koras-Russell cubic threefold} is the subvariety $X={\rm Spec}\left(A\right)$
of the affine space $\mathbb{A}^{4}={\rm Spec}\left(\mathbb{C}\left[x,y,z,t\right]\right)$
defined by the equation $x+x^{2}y+z^{2}+t^{3}=0$. It first appeared
in the work of Koras and Russell \cite{KoRu1,KoRu2} on the linearization
problem for algebraic actions of the multiplicative group $\mathbb{G}_{m}$
on the affine space $\mathbb{A}^{3}$. The question at that time was
to decide if $X$ is algebraically isomorphic to $\mathbb{A}^{3}$
or not, and a positive answer would have led to an example of a non
linearizable $\mathbb{G}_{m}$-action on $\mathbb{A}^{3}$. One of
the difficulties is that when equipped with the euclidean topology
$X$ is diffeomorphic to the euclidean space $\mathbb{R}^{6}$(see
e.g. \cite{ChouDim}). So it is impossible to distinguish $X$ from
$\mathbb{A}^{3}$ by topological invariants. Actually, it turned out
that all classical algebraic invariants fail to distinguish $X$ from
$\mathbb{A}^{3}$. 

Nowadays, the fact that $X$ is not algebraically isomorphic to $\mathbb{A}^{3}$
can be derived from a result of Kaliman \cite{Ka02} which says that
if the general fibers of regular function $f:\mathbb{A}^{3}\rightarrow\mathbb{A}^{1}$
are isomorphic to the affine plane $\mathbb{A}^{2}$, then \emph{all}
the closed fibers of $f$ are isomorphic to $\mathbb{A}^{2}$. On
the other hand, it is easily seen that the closed fibers of the projection
${\rm pr}_{x}:X\rightarrow\mathbb{A}^{1}$ are isomorphic to $\mathbb{A}^{2}$
except for ${\rm pr}_{x}^{-1}\left(0\right)$ which is isomorphic
to the cylinder $C\times\mathbb{A}^{1}$ over the cuspidal cubic curve
$C\simeq{\rm Spec}\left(\mathbb{C}\left[z,t\right]/\left(z^{2}+t^{3}\right)\right)$. 

But the problem was originally solved by Makar-Limanov \cite{ML96}
by a different method, based on the study of algebraic actions of
the additive group $\mathbb{G}_{a}$ on $X$. He established that
 $X$ is not algebraically isomorphic to $\mathbb{A}^{3}$ because
it admits ``fewer'' algebraic $\mathbb{G}_{a}$-actions than $\mathbb{A}^{3}$.
More precisely, Makar-Limanov introduced a new invariant of affine
algebraic varieties $V$ defined as the sub-algebra ${\rm ML}\left(V\right)$
of the coordinate ring of $V$ consisting of regular functions on
$V$ which are invariants under \emph{all} algebraic $\mathbb{G}_{a}$-actions
on $V$. For affine spaces, this invariant consists of constants only.
In contrast, Makar-Limanov established that ${\rm ML}\left(X\right)$
is isomorphic to the polynomial ring $\mathbb{C}\left[x\right]$.
To compute ${\rm ML}\left(X\right)$, Makar-Limanov used the correspondence
between algebraic $\mathbb{G}_{a}$-actions on an affine variety $V$
and \emph{locally nilpotent} $\mathbb{C}$-\emph{derivations} of its
coordinate ring $\mathbb{C}\left[V\right]$, that is, derivations
$\partial:\mathbb{C}\left[V\right]\rightarrow\mathbb{C}\left[V\right]$
such that every element of $\mathbb{C}\left[V\right]$ is annihilated
by a suitable power of $\partial$. Under this correspondence, $\mathbb{G}_{a}$-invariant
regular functions coincide with the elements of the kernel ${\rm Ker}\partial$
of the associated locally nilpotent derivation, and ${\rm ML}\left(V\right)$
can be equivalently defined as the intersection in $\mathbb{C}\left[V\right]$
of the kernels of all locally nilpotent derivations of $\mathbb{C}\left[V\right]$. 

It is easy to see that ${\rm ML}\left(X\right)\subset\mathbb{C}\left[x\right]$.
For instance, the locally nilpotent derivations $x^{2}\partial_{z}-2z\partial_{y}$
and $x^{2}\partial_{t}-3t^{2}\partial_{y}$ of $\mathbb{C}\left[x,y,z,t\right]$
annihilate the defining equation $x+x^{2}y+z^{2}+t^{3}=0$ of $X$
and induce non trivial locally nilpotent derivations $\partial_{1}$
and $\partial_{2}$ of the coordinate ring $A$ of $X$ such that
${\rm Ker}\left(\partial_{1}\right)\cap{\rm Ker}\left(\partial_{2}\right)=\mathbb{C}\left[x\right]$.
The main achievement of Makar-Limanov was to show that $\partial\left(x\right)=0$
for every locally nilpotent derivation of $A$. The original proof
has been simplified and generalized by many authors, but the key arguments
remain quite elaborate and depend on techniques of equivariant deformations
to reduce the problem to the study of homogeneous $\mathbb{G}_{a}$-actions
on certain affine cones associated with $X$ (see e.g., \cite{KML97},
\cite{KML07} and \cite{Zai99} ). 

Now, given a new variable $w$, we can identify $X\times\mathbb{A}^{1}$
with the subvariety of $\mathbb{A}^{5}={\rm Spec}\left(\mathbb{C}\left[x,y,z,t,w\right]\right)$
defined by the equation $x+x^{2}y+z^{2}+t^{3}=0$. Again, it is not
difficult to see that ${\rm ML}\left(X\times\mathbb{A}^{1}\right)\subset\mathbb{C}\left[x\right]$,
and it is natural to ask if ${\rm ML}\left(X\times\mathbb{A}^{1}\right)\not=\mathbb{C}$
or not. In turned out that Makar-Limanov techniques are inefficient
in this context, and very few progress has been made on this particular
problem since the late nineties. In this note, we prove the following
result.

\begin{thm*}
${\rm ML}\left(X\times\mathbb{A}^{1}\right)=\mathbb{C}$. 
\end{thm*}
\noindent A consequence of this result is that the Makar-Limanov
invariant carries no useful information to decide if $X\times\mathbb{A}^{1}$
is an \emph{exotic} $\mathbb{A}^{4}$, i.e. a variety diffeomorphic
to $\mathbb{R}^{8}$ but not algebraically isomorphic to $\mathbb{A}^{4}$. 

\begin{acknowledgement*}
The author is very grateful to Pierre-Marie Poloni for sharing his
key observation that the Koras-Russell cubic threefold can be interpreted
as a one-parameter family of Danielewski hypersurfaces and to Lucy
Moser-Jauslin for her very careful reading of a preliminary version
of the present note. 
\end{acknowledgement*}

\section{A Danielewski trick proof that ${\rm ML}\left(X\times\mathbb{A}^{1}\right)=\mathbb{C}$ }

Before giving the proof, we find it enlightening to review Danielewski's
classical counter-example to the Zariski Cancellation Problem. Indeed,
it formally contains, in a simpler form, all the ingredients needed
for the proof of the Theorem.

\subsection{Danielewski's construction }

\begin{enavant} Danielewski \cite{Dan89} established that the smooth
affine surfaces $S_{1}=\left\{ xz=y^{2}-1\right\} $ and $S_{2}=\left\{ x^{2}z=y^{2}-1\right\} $
in $\mathbb{A}^{3}={\rm Spec}\left(\mathbb{C}\left[x,y,z\right]\right)$
provide a counter example to the generalized Cancellation Problem,
that is, $S_{1}\times\mathbb{A}^{1}$ is isomorphic to $S_{2}\times\mathbb{A}^{1}$
but $S_{1}$ is not isomorphic to $S_{2}$. To show that $S_{1}\times\mathbb{A}^{1}$
is isomorphic to $S_{2}\times\mathbb{A}^{1}$, he exploited the fact
that $S_{1}$ and $S_{2}$ can be equipped with set-theoretically
free $\mathbb{G}_{a}$-actions induced by the $\mathbb{G}_{a}$-actions
on $\mathbb{A}^{3}$ associated with the locally nilpotent $\mathbb{C}\left[x\right]$-derivations
$x\partial_{y}+2y\partial_{z}$ and $x^{2}\partial_{y}+2y\partial_{z}$
of $\mathbb{C}\left[x,y,z\right]$ respectively. The fibers of the
$\mathbb{G}_{a}$-invariant projections ${\rm \pi}_{i}={\rm pr}_{x}\mid_{S_{i}}:S_{i}\rightarrow\mathbb{A}^{1}={\rm Spec}\left(\mathbb{C}\left[x\right]\right)$,
$i=1,2$ coincide with the orbits of the $\mathbb{G}_{a}$-actions
except $\pi_{i}^{-1}\left(0\right)$ which consists of the disjoint
union of two distinct orbits. In particular, $\pi_{i}:S_{i}\rightarrow\mathbb{A}^{1}$
is not a $\mathbb{G}_{a}$-bundle. However, Danielewski observed that
the $\pi_{i}$'s factor through Zariski locally trivial $\mathbb{G}_{a}$-bundles
$\rho_{i}:S_{i}\rightarrow\tilde{\mathbb{A}}^{1}$, $i=1,2$, over
the affine line with a double origin, obtained from $\mathbb{A}^{1}={\rm Spec}\left(\mathbb{C}\left[x\right]\right)$
by replacing its origin by two closed points, one for each of the
connected components of $\pi_{i}^{-1}\left(0\right)$. 

\end{enavant}

\begin{enavant} In turn, this implies that there exists a cartesian
diagram \[\xymatrix{ & S_1\times_{\tilde{A}^1} S_2 \ar[dl]_{{\rm pr}_1} \ar[dr]^{{\rm pr}_2} \\ S_1 \ar[dr]^{\rho_1} & & S_2 \ar[dl]_{\rho_2} \\ & \tilde{A}^1  }\]
where $S_{1}\times_{\tilde{\mathbb{A}}^{1}}S_{2}$ is a $\mathbb{G}_{a}$-bundle
over $S_{1}$ and $S_{2}$ via the first and the second projections
respectively. Since $S_{1}$ and $S_{2}$ are both affine, it follows
that $S_{1}\times_{\tilde{\mathbb{A}}^{1}}S_{2}$ is simultaneously
isomorphic to the trivial $\mathbb{G}_{a}$-bundles $S_{1}\times\mathbb{A}^{1}$
and $S_{2}\times\mathbb{A}^{1}$ over $S_{1}$ and $S_{2}$ respectively
(see e.g., XI.5.3 in \cite{SGA1}). This implies the existence of
an isomorphism $\Theta:S_{1}\times\mathbb{A}^{1}\stackrel{\sim}{\rightarrow}S_{2}\times\mathbb{A}^{1}$
of $\mathbb{A}^{2}$-bundles over $\tilde{\mathbb{A}}^{1}$, whence
of schemes over $\mathbb{A}^{1}\simeq{\rm Spec}\left(\Gamma(\tilde{\mathbb{A}}^{1},\mathcal{O}_{\tilde{\mathbb{A}}^{1}})\right)={\rm Spec}\left(\mathbb{C}\left[x\right]\right)$.

\end{enavant}

\begin{enavant} Although Danielewski argument was different, the
fact that $S_{2}$ and $S_{1}$ are not isomorphic can be deduced
from a result of Makar-Limanov \cite{ML01} asserting that ${\rm ML}\left(S_{2}\right)=\mathbb{C}\left[x\right]$,
together with the observation that due to the symmetry between the
variables $x$ and $z$ in the defining equation of $S_{1}$, one
has ${\rm ML}\left(S_{1}\right)=\mathbb{C}$. Since $S_{2}\times\mathbb{A}^{1}$
is isomorphic to $S_{1}\times\mathbb{A}^{1}$, it follows in particular
that ${\rm ML}\left(S_{2}\times\mathbb{A}^{1}\right)\simeq{\rm ML}\left(S_{1}\times\mathbb{A}^{1}\right)=\mathbb{C}$.
This can be reinterpreted more explicitly as follows. Certainly, one
has ${\rm ML}\left(S_{2}\times\mathbb{A}^{1}\right)\subset\mathbb{C}\left[x\right]$.
On the other hand, the locally nilpotent derivation $z\partial_{y}+2y\partial_{x}$
of $\mathbb{C}\left[x,y,z,w\right]$ induces a locally nilpotent derivation
$\delta_{1}$ of the coordinate ring $\mathbb{C}\left[x,y,z,w\right]/\left(xz-y^{2}+1\right)$
of $S_{1}\times\mathbb{A}^{1}$ such that $\delta_{1}\left(x\right)\neq0$.
Since $\Theta:S_{1}\times\mathbb{A}^{1}\stackrel{\sim}{\rightarrow}S_{2}\times\mathbb{A}^{1}$
is an isomorphism of schemes over ${\rm Spec}\left(\mathbb{C}\left[x\right]\right)$,
it follows that $\left(\Theta^{*}\right)^{-1}\delta_{1}\Theta^{*}$
is a locally nilpotent derivation $\delta$ of the coordinate ring
of $S_{2}\times\mathbb{A}^{1}$ such that $\delta\left(x\right)\neq0$,
and so, ${\rm ML}\left(S_{2}\times\mathbb{A}^{1}\right)=\mathbb{C}$. 

\end{enavant}

\subsection{Proof of the  Theorem}

\begin{enavant} For our purpose, it is more convenient to rewrite
the defining equation of $X={\rm Spec}\left(A\right)$ as $x^{2}z=y^{2}+x-t^{3}$.
This corresponds to making the coordinate change $\left(x,y,z,t\right)\mapsto\left(-x,z,iy,t\right)$
on the ambient space $\mathbb{A}^{4}$. As observed in the introduction,
one has certainly ${\rm ML}\left(X\times\mathbb{A}^{1}\right)\subset\mathbb{C}\left[x\right]$.
So ${\rm ML}\left(X\times\mathbb{A}^{1}\right)=\mathbb{C}$ provided
that we can find a locally nilpotent derivation $\partial$ of the
coordinate ring $A\left[w\right]$ of $X\times\mathbb{A}^{1}$ such
that $\partial x\neq0$. We may even suppose that we are looking for
such a derivation with the additional property that $\partial t=0$.
With this hypothesis, we can further reduce the problem to finding
a locally nilpotent $\mathbb{C}\left[t,t^{-1}\right]$-derivation
$\delta$ of \[
A\left[w\right]\otimes_{\mathbb{C}\left[t\right]}\mathbb{C}\left[t,t^{-1}\right]\simeq\mathbb{C}\left[x,y,z,t^{\pm1}\right]\left[w\right]/\left(x^{2}z-y^{2}-x+t^{3}\right)\]
such that $\delta\left(x\right)\neq0$. Indeed, since $A\left[w\right]$
is a finitely generated algebra, for a suitably chosen $k\geq0$,
$t^{k}\delta$ will extend to a locally nilpotent derivation $\partial$
of $A\left[w\right]$ such that $\partial\left(x\right)\neq0$. 

\end{enavant}

\begin{enavant} We let $Y_{*}={\rm Spec}\left(\mathbb{C}\left[x,t^{\pm1}\right]\right)\simeq\mathbb{A}^{1}\times\mathbb{A}_{*}^{1}$
and we consider the affine varieties $X_{1}={\rm Spec}\left(B_{1}\right)$
and $X_{2}=X\setminus\left\{ t=0\right\} ={\rm Spec}\left(B_{2}\right)$
where 

\[
B_{1}=\mathbb{C}\left[x,y,z,t^{\pm1}\right]/\left(xz-y^{2}+t^{3}\right)\quad\textrm{and}\quad B_{2}=\mathbb{C}\left[x,y,z,t^{\pm1}\right]/\left(x^{2}z-y^{2}-x+t^{3}\right)\]
The locally nilpotent $\mathbb{C}\left[x,t^{\pm1}\right]$-derivations
$x\partial_{y}+2y\partial_{z}$ and $x^{2}\partial_{y}+2y\partial_{z}$
of $\mathbb{C}\left[x,y,z,t^{\pm1}\right]$ induce locally nilpotent
derivations of $B_{1}$ and $B_{2}$ respectively, defining set-theoretically
free $\mathbb{G}_{a}$-actions $m_{i}:\mathbb{G}_{a}\times X_{i}\rightarrow X_{i}$,
$i=1,2$. The $\mathbb{G}_{a}$-equivariant projections $\pi_{i}={\rm pr}_{x,t}\mid_{X_{i}}:X_{i}\rightarrow Y_{*}$
restrict to trivial $\mathbb{G}_{a}$-bundles over $Y_{*}\setminus\left\{ x=0\right\} $.
In constrast, the fibers of the $\pi_{i}$'s over every closed point
of the punctured line $\left\{ x=0\right\} \subset Y_{*}$ consist
of the disjoint union of two $\mathbb{G}_{a}$-orbits, and their fiber
over the point $\left(x\right)\in Y_{*}={\rm Spec}\left(\mathbb{C}\left[x,t^{\pm1}\right]\right)$
is isomorphic to the affine line over the Galois extension $\mathbb{C}\left(t\right)\left[y\right]/\left(y^{2}-t^{3}\right)\simeq\mathbb{C}\left(t\right)\left[\mu\right]/\left(\mu^{2}-t\right)$
of the residue field $\kappa\left(\left(x\right)\right)=\mathbb{C}\left(t\right)$.
Informally, this indicates that $\pi_{i}:X_{i}\rightarrow Y_{*}$
should factor through a $\mathbb{G}_{a}$-bundle $\rho_{i}:X_{i}\rightarrow\mathfrak{S}$,
$i=1,2$ over a geometric object $\mathfrak{S}$ obtained from $Y_{*}$
by replacing the point $\left(x\right)$, i.e., the punctured line
$\left\{ x=0\right\} \subset Y_{*}$, not by two disjoint copies of
itself as in Danielewski's construction, but rather by a nontrivial
\'etale double covering of itself. 

\end{enavant}

\begin{enavant} Clearly, an object $\mathfrak{S}$ with the required
property cannot exist in the category of schemes. However, one can
construct such an $\mathfrak{S}$ in the larger category of algebraic
spaces as follows. We let $Z_{*}={\rm Spec}\left(\mathbb{C}\left[x,\mu^{\pm1}\right]\right)$
and we let $\mathfrak{S}$ be the quotient of $Z_{*}$ by the \'etale
equivalence relation $\left(x,\mu\right)\sim\left(x,-\mu\right)$
if $x\neq0$. More formally, this means that $\mathfrak{S}=Z_{*}/R$
where $\left(s,t\right):R\rightarrow Z_{*}\times Z_{*}$ is the \'etale
equivalence relation defined by \[
\left(s,t\right):R=Z_{*}\sqcup Z_{*}\setminus\left\{ x=0\right\} \rightarrow Z_{*}\times Z_{*}\left(x,\mu\right)\mapsto\begin{cases}
\left(\left(x,\mu\right),\left(x,\mu\right)\right) & \textrm{if }\left(x,\mu\right)\in Z_{*}\subset R\\
\left(\left(x,\mu\right),\left(x,-\mu\right)\right) & \textrm{if }\left(x,\mu\right)\in Z_{*}\setminus\left\{ x=0\right\} \subset R\end{cases}\]
Now the $R$-invariant morphism $Z_{*}\rightarrow Y_{*}$, $\left(x,\mu\right)\mapsto\left(x,\mu^{2}\right)$
descends to a morphism $\psi:\mathfrak{S}\rightarrow Y_{*}$ restricting
to an isomorphism outside $\left\{ x=0\right\} $ and with fiber over
$\left(x\right)$ isomorphic to ${\rm Spec}\left(\mathbb{C}\left(t\right)\left[\mu\right]/\left(\mu^{2}-t\right)\right)$
as desired. 

\end{enavant}

\begin{rem}
An alternative construction of $\mathfrak{S}$ is the following :
First we let $W$ be the scheme obtained by gluing two copies $W_{\pm}$
of $Z_{*}={\rm Spec}\left(\mathbb{C}\left[x,\mu^{\pm1}\right]\right)$
by the identity outside the punctured line $\left\{ x=0\right\} $.
The group $\mathbb{Z}_{2}$ acts freely on $W$ by $W_{\pm}\ni\left(x,\mu\right)\mapsto\left(x,-\mu\right)\in W_{\mp}$,
and $\mathfrak{S}$ coincides with the quotient $W/\mathbb{Z}_{2}$
taken in the category of algebraic spaces. Note that this $\mathbb{Z}_{2}$-action
is properly discontinuous in the analytic topology on $W$, so that
$\mathfrak{S}$ equipped with the quotient analytic topology has the
structure of a locally separated analytic space. 
\end{rem}
\begin{enavant} Let us assume for a moment that we have factorizations
\[
\pi_{i}=\psi\circ\rho_{i}:X_{i}\stackrel{\rho_{i}}{\rightarrow}\mathfrak{S}\stackrel{\psi}{\rightarrow}Y_{*},\; i=1,2\]
where $\rho_{i}:X_{i}\rightarrow\mathfrak{S}$, $i=1,2$ is an \'etale
locally trivial $\mathbb{G}_{a}$-bundle. Then $X_{1}\times_{\mathfrak{S}}X_{2}$
is an \'etale locally trivial $\mathbb{G}_{a}$-bundle over $X_{1}$
and $X_{2}$ via the first and the second projection respectively.
Again, these bundles are both trivial as $X_{1}$ and $X_{2}$ are
affine, and we obtain isomorphisms $X_{1}\times\mathbb{A}^{1}\stackrel{\sim}{\rightarrow}X_{1}\times_{\mathfrak{S}}X_{2}\stackrel{\sim}{\leftarrow}X_{2}\times\mathbb{A}^{1}$.
The induced isomorphism $\Theta:X_{1}\times\mathbb{A}^{1}\stackrel{\sim}{\rightarrow}X_{2}\times\mathbb{A}^{1}$
is an isomorphism of \'etale locally trivial $\mathbb{A}^{2}$-bundles
over $\mathfrak{S}$, whence, in particular, of schemes over $Y_{*}$.
Now the locally nilpotent $\mathbb{C}\left[t^{\pm1}\right]$-derivation
$2y\partial_{x}+z\partial_{y}$ of $\mathbb{C}\left[x,y,z,t^{\pm1},w\right]$
induces a locally nilpotent derivation $d$ of the coordinate ring
$B_{1}\left[w\right]$ of $X_{1}\times\mathbb{A}^{1}$ such that $d\left(x\right)\neq0$.
Since $\Theta:X_{1}\times\mathbb{A}^{1}\stackrel{\sim}{\rightarrow}X_{2}\times\mathbb{A}^{1}$
is an isomorphism of schemes over $Y_{*}={\rm Spec}\left(\mathbb{C}\left[x,t^{\pm1}\right]\right)$,
it follows that $\delta=\left(\Theta^{*}\right)^{-1}d\Theta^{*}$
is a locally nilpotent $\mathbb{C}\left[t^{\pm1}\right]$-derivation
of the coordinate ring $B_{2}\left[w\right]$ of $X_{2}\times\mathbb{A}^{1}$
such that $\delta\left(x\right)\neq0$. Combined with the previous
discussion, this shows that ${\rm ML}\left(X\times\mathbb{A}^{1}\right)=\mathbb{C}$. 

\end{enavant}

\begin{enavant} So it remains to check that the $\mathbb{G}_{a}$-invariant
morphisms $\pi_{i}:X_{i}\rightarrow Y_{*}$, $i=1,2$, admit the required
factorization. It is a standard fact that a set-theoretically free
$\mathbb{G}_{a}$-action on a scheme $V$ admits a categorical quotient
in the form of a $\mathbb{G}_{a}$-bundle $\rho:V\rightarrow V/\mathbb{G}_{a}$
over an algebraic space $V/\mathbb{G}_{a}$ (see e.g. 10.4 in \cite{LMB}).
Thus we only need to check that $X_{i}/\mathbb{G}_{a}\simeq\mathfrak{S}$,
$i=1,2$, and that the morphisms $\bar{\pi}_{i}:X_{i}/\mathbb{G}_{a}\rightarrow Y_{*}$
induced by the $\mathbb{G}_{a}$-invariant morphisms $\pi_{i}:X_{i}\rightarrow Y_{*}$
coincide with $\psi:\mathfrak{S}\rightarrow Y_{*}$. This can be seen
as follows. Letting $U=\mathbb{G}_{a}\times Z_{*}={\rm Spec}\left(\mathbb{C}\left[v\right]\left[x,\mu^{\pm1}\right]\right)$,
one checks first that the $\mathbb{G}_{a}$-equivariant morphisms
\begin{align*}
\phi_{1}:\mathbb{G}_{a}\times Z_{*}\rightarrow X_{1}, & \left(v,x,\mu\right)\mapsto\left(x,\mu^{3}+xv,2\mu^{3}v+xv^{2},\mu^{2}\right)\\
\phi_{2}:\mathbb{G}_{a}\times Z_{*}\rightarrow X_{2}, & \left(v,x,\mu\right)\mapsto\left(x,\mu^{3}-\frac{1}{2\mu^{3}}x+x^{2}v,\frac{1}{4\mu^{6}}+\left(2\mu^{3}-\mu^{-3}x\right)v+x^{2}v^{2},\mu^{2}\right)\end{align*}
define \'etale trivializations of the $\mathbb{G}_{a}$-actions $m_{i}:\mathbb{G}_{a}\times X_{i}\rightarrow X_{i}$
on $X_{i}$, $i=1,2$. Then one checks easily that we have $\mathbb{G}_{a}$-equivariant
isomorphisms\begin{align*}
\xi_{1}:\mathbb{G}_{a}\times R & \stackrel{\sim}{\longrightarrow}U\times_{X_{1}}U\\
\left(x,\mu\right)\mapsto & \begin{cases}
\left(\phi_{1}\left(v,x,\mu\right),\phi_{1}\left(v,x,\mu\right)\right) & \textrm{if }\left(v,x,\mu\right)\in\mathbb{G}_{a}\times Z_{*}\subset R\\
\left(\phi_{1}\left(v,x,\mu\right),\phi_{1}\left(v+2\mu^{3}x^{-1},x,-\mu\right)\right) & \textrm{if }\left(v,x,\mu\right)\in\mathbb{G}_{a}\times Z_{*}\setminus\left\{ x=0\right\} \subset R\end{cases}\end{align*}
and \begin{align*}
\xi_{2}:\mathbb{G}_{a}\times R & \stackrel{\sim}{\longrightarrow}U\times_{X_{2}}U\\
\left(v,x,\mu\right) & \mapsto\begin{cases}
\left(\phi_{2}\left(v,x,\mu\right),\phi_{2}\left(v,x,\mu\right)\right) & \textrm{if }\left(v,x,\mu\right)\in\mathbb{G}_{a}\times Z_{*}\subset R\\
\left(\phi_{2}\left(v,x,\mu\right),\phi_{2}\left(v-\mu^{-3}x^{-1}+2\mu^{3}x^{-2},x,-\mu\right)\right) & \textrm{if }\left(v,x,\mu\right)\in\mathbb{G}_{a}\times Z_{*}\setminus\left\{ x=0\right\} \subset R.\end{cases}\end{align*}
By construction, the projections $\left({\rm pr}_{1},{\rm pr}_{2}\right):U\times_{X_{i}}U\rightrightarrows U=\mathbb{G}_{a}\times Z_{*}$
are \'etale and descend to the ones $\left(s,t\right):R\simeq U\times_{X_{i}}U/\mathbb{G}_{a}\rightrightarrows Z_{*}=\mathbb{G}_{a}\times Z_{*}/\mathbb{G}_{a}$
in such a way that we have a cartesian diagram \[\xymatrix{U\times_{X_i} U \ar@<0.5ex>[r]^-{{\rm pr}_1} \ar@<-0.5ex>[r]_-{{\rm pr}_2} \ar[d]& U \ar[d] \\ R=U\times_{X_i} U/\mathbb{G}_a \ar@<0.5ex>[r]^{s} \ar@<-0.5ex>[r]_{t} & Z_*=U/\mathbb{G}_a }\] 
Since $X_{i}$ coincides with the quotient of $U\times_{X_{i}}U$
by the \'etale equivalence relation $\left({\rm pr}_{1},{\rm pr}_{2}\right):U\times_{X_{i}}U\rightrightarrows U$,
it follows from I.5.8 in \cite{Kn} that the $\mathbb{G}_{a}$-bundle
$U\rightarrow Z_{*}=U/\mathbb{G}_{a}$ descends to a morphism of algebraic
spaces $\rho_{i}:X_{i}\rightarrow\mathfrak{S}=Z_{*}/R$, and that
we have a commutative diagram \[\xymatrix{U\times_{X_i} U \ar@<0.5ex>[r]^-{{\rm pr}_1} \ar@<-0.5ex>[r]_-{{\rm pr}_2} \ar[d]& U \ar[d] \ar[r] & X_i \ar[d]^{\rho_i}\\ R \ar@<0.5ex>[r]^-{s} \ar@<-0.5ex>[r]_-{t} & Z_*=U/\mathbb{G}_a \ar[r] & \mathfrak{S}=Z_*/R}\] 
in which the right hand side square is cartesian. This implies that
$\rho_{i}:X_{i}\rightarrow\mathfrak{S}$ is an \'etale locally trivial
$\mathbb{G}_{a}$-bundle, which shows that $\mathfrak{S}$ is isomorphic
to $X_{i}/\mathbb{G}_{a}$, $i=1,2$ as desired. Now the fact that
$\pi_{i}:X_{i}\rightarrow Y_{*}$ factors as $\psi\circ\rho_{i}$,
$i=1,2$, follows trivially from the construction. 

\end{enavant}

\begin{rem}
The maps $\rho_{i}:X_{i}\rightarrow\mathfrak{S}$, $i=1,2$, are holomorphic
$\mathbb{G}_{a}$-bundles when the $X_{i}$'s and $\mathfrak{S}$
are equipped with the analytic topology. Indeed, one can check that
the $\mathbb{G}_{a}$-invariant maps ${\rm pr}_{2}\mid_{\tilde{X}_{i}}:\tilde{X}_{i}=X_{i}\times_{Y_{*}}Z_{*}\rightarrow Z_{*}$
obtained from the base change by the \'etale Galois covering $Z_{*}\rightarrow Y_{*}$,
$\left(x,\mu\right)\mapsto\left(x,\mu^{2}\right)$ factor through
$\mathbb{Z}_{2}$-equivariant holomorphic $\mathbb{G}_{a}$-bundles
$\tilde{\rho}_{i}:\tilde{X}_{i}\rightarrow W$ such that $\rho_{i}=\tilde{\rho}_{i}/\mathbb{Z}_{2}:X_{i}\simeq\tilde{X}_{i}/\mathbb{Z}_{2}\rightarrow W/\mathbb{Z}_{2}\simeq\mathfrak{S}$,
$i=1,2$. 
\end{rem}
\noindent

\begin{rem}
The above descriptions imply that the isomorphy classes of the $\mathbb{G}_{a}$-bundles
$\rho_{1}:X_{1}\rightarrow\mathfrak{S}$ and $\rho_{2}:X_{2}\rightarrow\mathfrak{S}$
in $H_{\textrm{ét}}^{1}\left(\mathfrak{S},\mathbb{G}_{a}\right)\simeq H_{\textrm{ét}}^{1}\left(\mathfrak{S},\mathcal{O}_{\mathfrak{S}}\right)$
are represented by the non cohomologous \v{C}ech $1$-cocycles \[
\left\{ 0,2\mu^{3}x^{-1}\right\} \in\Gamma\left(R,\mathcal{O}_{R}\right)\quad\textrm{and}\quad\left\{ 0,-\mu^{-3}x^{-1}+2\mu^{3}x^{-2}\right\} \in\Gamma\left(R,\mathcal{O}_{R}\right)\]
for the \'etale covering $Z_{*}\rightarrow\mathfrak{S}$ . So the
varieties $X_{1}$ and $X_{2}$ are not isomorphic as $\mathbb{G}_{a}$-bundles
over $\mathfrak{S}$. Actually, one can check that ${\rm ML}\left(X_{1}\right)=\mathbb{C}\left[t^{\pm1}\right]$
whereas ${\rm ML}\left(X_{2}\right)=\mathbb{C}\left[t^{\pm1}\right]\left[x\right]$,
so that $X_{1}$ and $X_{2}$ are not even isomorphic as abstract
affine varieties. Thus they provide a counter-example to the Cancellation
Problem for factorial affine threefolds (see \cite{FM08} for other
counter-examples).

\bibliographystyle{amsplain}

\end{rem}

\end{document}